\documentclass[a4paper,11pt]{article}
\usepackage[utf8x]{inputenc}
\usepackage{amsthm,amsmath,amssymb,oldgerm}
\usepackage[margin=3cm]{geometry}
\theoremstyle{plain}
\newtheorem{thm}{Theorem}

\newtheorem{cor}[thm]{Corollary}
\theoremstyle{definition}
 
\theoremstyle{definition}
\newtheorem{rem}[thm]{Remark}

\DeclareMathOperator{\bk} {{\it{\mathcal{B}_{X}^{\tilde{k}}}}}
\DeclareMathOperator{\bkv} {{\it{\mathcal{B}_{X,\nu}^{\tilde{k}}}}}
\DeclareMathOperator{\bkappav} {{\it{\mathcal{B}_{X,\nu}^{\tilde{\kappa}}}}}

\DeclareMathOperator{\bkm} {{\it{\mathcal{B}_{M,\mathcal{L}}^{k}}}}
\DeclareMathOperator{\hkm} {{\it{K_{M,\mathcal{L}}^{k}}}}

\DeclareMathOperator{\hyp}{\mu_{\mathbb{H}}} 
\DeclareMathOperator{\hypnvol}{\mu_{\mathbb{H}}^{vol}} 
\DeclareMathOperator{\hypn}{\mu_{\mathbb{H}^{{\it{n}}}}} 
\DeclareMathOperator{\shyp}{\mu_{shyp}}

\DeclareMathOperator{\sk}{{\it{S^{\tilde{k}}(\Gamma)}}}
\DeclareMathOperator{\sknu}{{\it{S^{\tilde{k}}_{\nu}(\Gamma)}}}

\DeclareMathOperator{\tk}{{\it{\tilde{k}}}}
\let\k\relax
\DeclareMathOperator{\k}{{\it{(k,\ldots,k)}}}

\DeclareMathOperator{\lk}{{\it{\overline{\square}_{k}}}}
\let\ell\relax
\DeclareMathOperator{\ell}{\mathcal{L}}
\DeclareMathOperator{\hn}{\mathbb{H}^{{\it{n}}}}
\DeclareMathOperator{\zn}{{\it{(z_1,\ldots,z_n)}}}
\title{Estimates of Hilbert modular cusp forms of half-integral and integral weight}
{\small\author{Anilatmaja Aryasomayajula}}
\date{}
\begin{document}
\maketitle
\begin{abstract}\noindent  Let $\Gamma$ be a cocompact, discrete, and irreducible subgroup of 
$\mathrm{PSL}_{2}(\mathbb{R})^{n}$. Let $\nu$ be a unitary character of $\Gamma$. For 
$k\in1\slash 2\,\mathbb{Z}$, let $\sknu$ denote the complex vector space of cusp forms of weight-$\tk=\k$ and 
nebentypus $\nu^{2k}$ with respect to $\Gamma$. We assume that $\omega_{X,\nu}$, the line bundle 
of cusp forms of weight-$\tilde{1\slash 2}:=(1\slash 2,\ldots,1\slash2)$ with nebentypus $\nu$ over $X$ 
exists. Let $\lbrace f_{1},\ldots,f_{j_{\tk}} \rbrace$ denote an orthonormal basis of $\sknu$. In this article, 
we show that as $k\rightarrow \infty$, the sum $\sum_{i=1}^{j_{\tk}}y^{k}|f_{i}(z)|^{2}$ is bounded by 
$O(k^{n})$, where  the implied constant is independent of $\Gamma$. Furthermore, we extend these results to 
the case when $k\in2\mathbb{Z}$, and to the case when $\Gamma$ is commensurable with the Hilbert modular 
group $\Gamma_{K}:=\mathrm{PSL}_{2}(O_{K})$, where $K$ is a totally real number field of degree $n\geq 2$, 
and $\mathcal{O}_{K}$ is the ring of integers of $K$.   

\vspace{0.2cm}\noindent
\end{abstract}
\section{Introduction}\label{introduction}
The main aim of this article is to extend the heat kernel approach from \cite{anil} to study Hilbert modular 
cusp forms. J.~Jorgenson and J.~Kramer have successfully applied heat kernel analysis to study automorphic 
forms, and derived optimal estimates for various arithmetic invariants. In \cite{anil}, using a different 
approach from geometric analysis, we derived similar estimates for integral weight cusp forms 
as in \cite{jk2}, and extended these estimates to half-integral weight cusp forms. 

\vspace{0.2cm}\noindent
In the current article we extend the methods from \cite{anil} to study estimates of half-integral and integral 
weight Hilbert modular forms. The method can be extended to study automorphic forms of arbitrary Shimura varieties.  

\vspace{0.025cm}
\paragraph{Notation}
Let $\mathbb{H}$ denote the hyperbolic upper half-plane. Let $\hyp$ denote the $(1,1)$-form on $\mathbb{H}$ 
corresponding to the hyperbolic metric. For $z=x+iy\in\mathbb{H}$, it is given by the following formula
\begin{equation*}
 \hyp(z):= \frac{i}{2}\cdot\frac{dz\wedge d\overline{z}}{y^{2}}.
\end{equation*} 
For $\zn\in\hn$, let $\hypn$ denote the $(1,1)$-form on $\mathbb{H}^{n}$, given by the following formula
\begin{align*}
 \hypn(z):=\sum_{i=1}^{n}\hyp(z_i).
\end{align*}
Let $\hypnvol$ denote volume form corresponding to $\hypn$. 

\vspace{0.2cm}\noindent
Let $\Gamma\subset \mathrm{PSL}_{2}(\mathbb{R})^{n}$ be a cocompact, discrete, and irreducible (as in the sense 
of \cite{freitag}, and described in section \ref{section3}) subgroup acting on $\mathbb{H}^{n}$, via 
fractional linear transformations. We assume that $\Gamma$ admits no elliptic elements. So, the quotient space $X:=\Gamma\backslash 
\mathbb{H}$ admits the structure of a compact complex manifold of complex dimension $n$. We further assume that $\hypn$ 
defines a K\"ahler metric on $X$, which is compatible with the complex structure of $X$. 

\vspace{0.2cm}
Let $\nu$ denote a unitary character associated to $\Gamma$. For $k\in 1\slash 2\,\mathbb{Z}$, let $\sknu$ 
denote the complex vector space of cusp forms of weight-$\tk=\k$ and nebentypus 
$\nu^{2k}$ with respect to $\Gamma$. Similarly, for $k\in 2\mathbb{Z}$, let $\sknu$ denote the complex vector 
space of cusp forms of weight-$\tk=\k$ and nebentypus $\nu^{k\slash 2}$ with respect to $\Gamma$. 

For $k\in 1\slash 2\,\mathbb{Z}$ (or $k\in2\mathbb{Z}$), when the nebentypus $\nu$ is trivial, we denote $\sknu$ by 
$\sk$. 

\vspace{0.2cm}
For $k\in 1\slash 2\,\mathbb{Z}$ (or $k\in2\mathbb{Z}$), let $j_{\tk}$ denote the dimension of $\sknu$, and let 
$\lbrace f_{1},\ldots,f_{j_{\tk}} \rbrace$ denote an orthonormal basis of $\sknu$ with respect to the Petersson 
inner-product. For any $z=(x_1+iy_1,\ldots,x_n+iy_n)\in\hn$, put 
\begin{align*}
\bkv(z):=\sum_{i=1}^{j_{\tk}}\bigg(\prod_{j=1}^{n}y_j^{k}|f_{i}(z)|^{2}\bigg).
\end{align*}
When the nebentypus $\nu$ is trivial, we denote $\bkv(z)$ by $\bk(z)$. 

\vspace{0.2cm}
We further assume that $\omega_{X,\nu}$ and $\Omega_{X,\nu}$, the line bundles of cusp forms of 
weight-$\tilde{2}:=(2,\ldots,2)$  and weight-$\tilde{1\slash 2}:=(1\slash2,\ldots, 1\slash 2)$, respectively, 
with nebentypus $\nu$, over $X$ exist. 

\vspace{0.2cm}
With notation as above, for $k\in1\slash 2\,\mathbb{Z}$, (or $k\in 2\mathbb{Z}$), we prove the following estimate
\begin{align}\label{estimate1}
\lim_{k\rightarrow \infty}\sup_{z\in X}\frac{1}{k^n}\bkv(z)= O(1),
\end{align}
where the implied constant is independent of $\Gamma$. 

\vspace{0.2cm}
Let $K$ be a totally real number field of degree $n\geq 2$, and let $\mathcal{O}_{K}$ denote its ring of integers. 
Let $\Gamma_{K}:=\mathrm{PSL}_{2}(\mathcal{O}_{K})$ denote the Hilbert modular group, and let $\Gamma$ now be 
commensurable with $\Gamma_{K}$. We further assume that $\Gamma$ does not admit elliptic elements, which implies 
that $X=\Gamma \backslash\hn$ admits the structure of a noncompact complex manifold of complex dimension $n$ 
with K\"ahler metric $\hypn$. Then, for any $z\in X$ and $k\in1\slash 2\,\mathbb{Z}$, we prove the following estimate
\begin{align}\label{estimate2}
\limsup_{k\rightarrow \infty}\frac{1}{k^{n}}\bkv(z)\leq\frac{1}{8\pi};
\end{align} 
and for $k\in 2\mathbb{Z}$, we prove that
\begin{align}\label{estimate3}
\limsup_{k\rightarrow \infty}\frac{1}{k^{n}}\bkv(z)\leq\frac{1}{2\pi}.
\end{align}  
\vspace{0.025cm}
\paragraph{Applications and existing results}
When $X$ is compact, in proving estimate \eqref{estimate1}, for any $z\in X$ and 
$k\in 1\slash 2\,\mathbb{Z}$, we first arrive at the following estimate
\begin{align*}
\lim_{k\rightarrow \infty}\frac{1}{k^{n}}\bkv(z) =O(1),
\end{align*}
where the implied constant is independent of $\Gamma$. Furthermore, for $k\gg0$, adapting the 
Selberg trace formula method from \cite{freitag} (Theorem 1.6 on p. 79) to compute $j_{\tk}$ the dimension 
of $\sknu$, one can show that $j_{\tk}=O(k^{n})$. Using which, we derive that 
\begin{align*}
\lim_{k\rightarrow \infty}\frac{1}{j_{k}}\bkv(z) =O(1).
\end{align*} 
So let
\begin{align*}
\lim_{k\rightarrow \infty}\frac{1}{j_k}\bkv(z)=C.
\end{align*}
Then, for any $z\in X$ and $k\in1\slash 2\,\mathbb{Z}$, we observe that
\begin{align*}
 \lim_{k\rightarrow \infty}\frac{1}{Cj_k}\bkv(z)\shyp(z)=\shyp(z),
\end{align*}
which proves an equidistribution result for a set of orthonormal basis of Hilbert modular cusp forms. 
\vspace{0.2cm}

Let $k\in 2\mathbb{Z}$ and $\Gamma$ be commensurable with the Hilbert modular group $\mathrm{PSL}_{2}(
\mathcal{O}_{K})$. Furthermore, let the nebentypus $\nu$ be trivial. Then, using an infinite series 
representation for $\bk(z)$, Liu has proved the above equidistribution result in \cite{liu}. Using 
Liu's technique, in \cite{codgell}, Codgell and Luo extended the equidistribution result to Siegel 
modular cusp forms. 
\vspace{0.2cm}

Sup-norm bounds for automorphic forms, especially, the ones associated to Fuchsian subgroups of 
$\mathrm{PSL}_{2}(\mathbb{R})$ is a  widely studied topic in number theory. However, sup-norm bounds 
for Hilbert modular cusp forms are yet to be explored. But in principle, most of the methods used to study 
sup-norm bounds for automorphic forms associated to Fuchsian subgroups should extend  to higher dimensions. 
\vspace{0.2cm}

We now state relevant results associated to Fuchsian subgroups of $\mathrm{PSL}_{2}(\mathbb{R})$. 
For the rest of this section, let $N\in\mathbb{N}$ with $N$ square-free.

Let $f$ any Hecke-normalized newform of $\Gamma_{0}(N)$ with trivial nebentypus and of weight $2$. 
Then, in \cite{abbes}, Abbes and Ullmo proved the following estimate 
\begin{align}\label{abbesestimate}
\sup_{z\in\mathbb{H}} y|f(z)|=O_{\varepsilon}(N^{\frac{1}{2}+\varepsilon}),
\end{align}
for any $\varepsilon > 0$, and the implied constant depends on $\varepsilon$. 
\vspace{0.2cm}

Let $f$ any Hecke-Maass cuspidal newform of $\Gamma_{0}(N)$ with trivial nebentypus and Laplacian 
eigenvalue $\lambda >0$. Then, in \cite{temp}, Harcos and Templier proved the 
following estimate 
\begin{align*}
\sup_{z\in\mathbb{H}} |f(z)|=O_{\lambda,\varepsilon}(N^{-\frac{1}{6}+\varepsilon}),
\end{align*}
for any $\varepsilon > 0$, and the implied constant depends on the eigenvalue $\lambda$ and $\varepsilon$. 

\vspace{0.2cm}
For $k\in1\slash 2\,\mathbb{Z}$, let $f$ be any cusp form of $\Gamma_{0}(4N)$ with nebentypus $\nu$ 
and of weight-$k$. Furthermore, let $f$ be normalized with respect to the Petersson inner-product. Then, in 
\cite{kiral}, Kiral has derived the following estimate 
\begin{align*}
 \sup_{z\in Y_{0}(N)}y^{k}|f(z)|^{2}=O_{k,\varepsilon}\big(N^{\frac{1}{2}-\frac{1}{18}+\varepsilon}\big),
\end{align*}
for any $\varepsilon>0$, and the implied constant depends on the weight $k$ and $\varepsilon$. 

\vspace{0.2cm}
Using heat kernel techniques, Jorgenson, and Kramer have re-proved the result of Abbes and Ullmo 
(estimate \eqref{abbesestimate} in \cite{jk1}). In \cite{jk2}, Friedman, Jorgenson and Kramer 
extended their method from \cite{jk1}, to derive sup-norm bounds for $\bk(z)$. When $X$ is a compact 
hyperbolic Riemann surface, they proved the following estimate
\begin{align*}
 \sup_{z\in X}\bk(z)=O(k),
\end{align*}
where the implied constant is independent of $X$. When $X$ is a noncompact hyperbolic Riemann 
surface of finite volume, they proved the following estimate
\begin{align}\label{jkestimate}
 \sup_{z\in X}\bk(z)=O(k^{\frac{3}{2}}),
\end{align}
where the implied constant is independent of $X$. The estimates of Jorgenson and Kramer are 
optimal, as shown in \cite{jk2}. 

\vspace{0.2cm}
For $k\in\mathbb{R}_{>0}$ with $k>2$, $\bkv(z)$ can be represented by an infinite series, 
which is uniformly convergent in $z\in X$. Using which, Steiner has extended the bounds of Jorgenson and 
Kramer to real weights and nontrivial nebentypus. 

\vspace{0.2cm} 
Let $\Gamma$ be any subgroup of finite index in $\mathrm{SL}_{2}(\mathbb{Z})$, and let $A$ be a compact subset of $X$. 
Then, in \cite{raphael}, for $k\in\mathbb{R}_{>0}$ with $k\gg1$, Steiner has derived the following estimates
\begin{align}\label{sestimate}
 \sup_{z\in A}\bkv(z)=O_{A}(k),
\end{align}
where the implied constant depends on the compact subset $A$; and 
\begin{align*}
\sup_{z\in X }\bkv(z)=O_{X}\big(k^{\frac{3}{2}}\big),
\end{align*}
where the implied constant depends on $X$. 

\vspace{0.2cm}
It is difficult to directly extend the methods of Ullmo and Abbes, Harcos and Templier, and that of Kiral 
to Hilbert modular cusp forms. It is possible to extend Steiner's method, but the results could be 
restricted to groups commensurable with the Hilbert modular group.  

\vspace{0.2cm}
The method of Jorgenson and Kramer relies on the analysis of the infinite series representation of the 
heat kernel associated to the Riemann surface. It is possible to extend their heat kernel analysis to higher 
dimensions, namely to study Hilbert modular cusp forms and Siegel modular cusp forms of integral weight. 
However, one has to address certain nontrivial convergence issues, while doing so. Their method cannot be 
extended to cusp forms with nontrivial nebentypus, nor to the case of cusp forms of half-integral weight. 
\section{Heat kernels on compact complex manifolds}
In this section, we recall the main results from \cite{bouche} and \cite{berman}, which we use in 
the next section. 

\vspace{0.2cm}
Let $(M,\omega)$ be a compact complex manifold of dimension $n$ with natural Hermitian metric $\omega$. Let 
$\ell$ be a positive Hermitian holomorphic line bundle on $M$ with the Hermitian metric given by 
$\|s(z)\|^{2}_{\ell}:=e^{-\phi(z)}|s(z)|^{2}$, where $s\in\ell$ is any section, and $\phi(z)$ is 
a real-valued function defined on $M$. 

\vspace{0.2cm}
For any $k\in\mathbb{N}$, let $\lk:=(\overline{\partial}^{\ast}+\overline{\partial})^{2}$ denote the 
$\overline{\partial}$-Laplacian acting on smooth sections of the line bundle $\ell^{\otimes k}$. Let $\hkm(t;z,w)$ denote the smooth kernel of the 
operator $e^{-\frac{2t}{k}\lk}$. We refer the reader to p. 2 in \cite{bouche}, for the details regarding the 
properties which uniquely characterize the heat kernel $\hkm(t;z,w)$. When $z=w \in M$, the heat kernel 
$\hkm(t;z,z)$ admits the following spectral expansion
\begin{align}\label{spectralexpn}
\hkm(t;z,z)=\sum_{n\geq 0} e^{-\frac{2t}{k}\lambda_{n}^{k}}\|\varphi_{n}(z)\|_{\mathcal{L}^{\otimes k}}^{2},
\end{align}
where $\lbrace\lambda_{n}^k\rbrace_{n\geq0}$ denotes the set of eigenvalues of $\lk$ 
(counted with multiplicities), and $\lbrace\varphi_{n}\rbrace_{n\geq0}$ denotes a set of associated 
orthonormal eigenfunctions. 

\vspace{0.2cm}
Let $H^{0}(M,\mathcal{L}^{\otimes k})$ denote the vector space of global holomorphic sections of the line bundle 
$\mathcal{L}^{\otimes k}$, and let  $\lbrace s_{i}\rbrace$ denote an orthonormal basis of 
$H^{0}(M,\mathcal{L}^{\otimes k})$. For any $z\in M$, the following function is called the Bergman kernel 
associated to the line bundle $\mathcal{L}^{\otimes k}$
\begin{align}\label{bkdefn}
\bkm(z):= \sum_{i}\| s_{i}(z)\|_{\mathcal{L}^{\otimes k}}^{2}.
\end{align}
The above definition is independent of the choice of orthonormal basis of $H^{0}(M,\mathcal{L}^{\otimes k})$. 

\vspace{0.2cm}
For any $z\in M$ and $t\in\mathbb{R}_{>0}$, from the spectral expansion of the heat kernel $\hkm(t;z,w)$ 
described in equation \eqref{spectralexpn}, it is easy to see that
\begin{align}\label{hkbkreln}
\bkm(z)\leq \hkm(t;z,z)\quad \mathrm{and}\quad
\lim_{t\rightarrow\infty} \hkm(t;z,z)=\bkm(z).
\end{align}
Let 
\begin{align}\label{curvatureform}
c_{1}(\ell)(z):=\frac{i}{2\pi}\partial\overline{\partial}\phi(z)
\end{align}
denote the curvature form of the line bundle $\ell$ at the point $z\in M$. Let $\alpha_{1},\ldots,\alpha_{n}$ 
denote the eigenvalues of $\partial\overline{\partial}\phi(z)$ at the point $z\in M$. Then, with notation 
as above, from Theorem 1.1 in \cite{bouche}, for any $z\in M$ and $t\in (0,k^{\varepsilon})$, for 
a given $\varepsilon>0$ not depending on $k$, we have
\begin{align}\label{boucheeqn1}
\lim_{k\rightarrow\infty}\frac{1}{k^{n}}\hkm(t;z,z)=\prod_{j=1}^{n}\frac{\alpha_{j}}{(4\pi)^{n}\sinh(\alpha_{j}t)}, 
\end{align}
and the convergence of the above limit is uniform in $z$. 

\vspace{0.2cm}
Using equations \eqref{hkbkreln} and \eqref{boucheeqn1}, in Theorem 2.1 in \cite{bouche}, Bouche derived the 
following estimate 

\vspace{0.025cm}
\begin{align}\label{boucheeqn2}
\lim_{k\rightarrow\infty}\frac{1}{k^{n}}\bkm(z)= O\big(\mathrm{det}_{\omega}\big(c_{1}(\ell)(z)\big)\big),
\end{align}
\vspace{0.01cm}

and the convergence of the above limit is uniform in $z\in M$. 

\vspace{0.2cm}
When $M$ is a noncompact complex manifold, using micro-local analysis of the Bergman kernel, in \cite{berman}, 
Berman derived the following estimate

\vspace{0.025cm}
\begin{align}\label{bermaneqn}
\limsup_{k\rightarrow\infty}\frac{1}{k^{n}}\bkm(z)\leq \mathrm{det}_{\omega}\big(c_{1}(\ell)(z)\big).
\end{align}
\vspace{0.01cm}

Unlike Bouche, Berman worked directly with the Bergman kernel $\bkm(z)$. His proof relies on comparison of 
$\bkm$ with the Bergman kernel of the trivial line bundle with constant metric defined over $\mathbb{C}^n$. 
He then used certain strong inequalities from the theory of elliptic partial differential equations to 
derive estimate \eqref{bermaneqn}. 
\section{Estimates of cusp forms}\label{section3}
As in section \ref{introduction}, let  $K$ be a totally real number field of degree $n\geq 2$, and let $\mathcal{O}_{K}$ be its 
ring of integers. Let $\sigma_{1},\ldots,\sigma_{n}$ be the $n$-real embeddings of the number field $K$. 
Then, the the Hilbert modular group $\Gamma_{K}:=\mathrm{PSL}_{2}(\mathcal{O}_{K})$ can be seen 
as a discrete subgroup of $\mathrm{PSL}_{2}(\mathbb{R})^{n}$, via the following map
\begin{align*}
\mathfrak{i}:\Gamma_{K}\hookrightarrow &\mathrm{PSL}_{2}(\mathbb{R})^{n}\\
\big(a_{ij}\big)_{1\leq i,j\leq n}\mapsto &\bigg(\big(\sigma_{1}(a_{ij})\big)_{1\leq i,j\leq n},\ldots,
\big(\sigma_{n}(a_{ij})\big)_{1\leq i,j\leq n} \bigg).
\end{align*}
We say that $\Gamma$ is commensurable with $\Gamma_{K}$, if there exists a $g\in\mathrm{PSL}_{2}
(\mathbb{R})^{n}$ such that $g\Gamma g^{-1}$ is a finite index subgroup of $\mathfrak{i}(\Gamma_{K})$. 

\vspace{0.2cm}
Any $\gamma:=(\gamma_{1},\ldots,\gamma_{n})\in \mathrm{PSL}_{2}(\mathbb{R})^{n}$, where 
$$
\gamma_{i}:=\bigg( 
\begin{array}{cc}
a_{i} & b_{i} \\
c_{i} &  d_{i} \\
\end{array}\bigg)\in \mathrm{PSL}_{2}(\mathbb{R})
$$
defines the following linear transformation on $\hn$
\begin{align*}
\gamma:\hn\longrightarrow &\hn\\
\zn\hookrightarrow&\big(\gamma_{1}(z_{1}),\ldots,\gamma_{n}(z_{n})\big),\,\,\mathrm{where}\,\,
\gamma_{i}(z_{i}):=\frac{a_{i}z_{i}+b_{i}}{c_{i}z_{i}+d_{i}}. 
\end{align*}
A subgroup $\Gamma\subset \mathrm{PSL}_{2}(\mathbb{R})^{n}$ is irreducible, if the restriction of 
each of the $n$ projections
\begin{align*}
\rho_{j}: \mathrm{PSL}_{2}(\mathbb{R})^{n}\longrightarrow \mathrm{PSL}_{2}(\mathbb{R}) \,\,(1\leq j\leq n)
\end{align*}
to $\Gamma$ is injective. 

\vspace{0.2cm}
Let $\Gamma$ be a discrete, irreducible subgroup of $\mathrm{PSL}_{2}(\mathbb{R})^{n}$, with no elliptic elements. 
Let $X:=\Gamma\backslash\hn$ denote the quotient space. For the rest of this article, we assume that 
$\Gamma$ is one of the following two types:

\vspace{0.1cm}
{\it{Type(1)}}: The group $\Gamma$ is cocompact, which implies that $X$ admits the structure of a compact complex 
manifold of complex dimension $n$. We assume that $\hypn$, the natural metric on $\hn$ defines a 
K\"ahler metric on $X$. 

\vspace{0.15cm}
{\it{Type(2)}}: The group $\Gamma$ is commensurable with the Hilbert modular group $\Gamma_{K}$, which implies that 
$X$ admits the structure of a noncompact complex manifold of complex dimension $n$ with cusps. The number 
of cusps is equal to the class number of $K$. The hyperbolic metric $\hypn$, which is the natural metric on $\hn$ defines a 
K\"ahler metric on $X$.

\vspace{0.2cm}
For the convenience of the reader, in this section, we adopt the following notation. We denote half-integers 
by $\kappa$, and even integers by $k$. 

\vspace{0.2cm}
Let $\nu$ be a unitary character of the group $\Gamma$. We assume that $\omega_{X,\,\nu}$, the line bundle of cusp 
forms of weight-$\tilde{1\slash 2}=(1\slash2,\ldots,1\slash 2)$ with nebentypus $\nu$ over $X$ exists. 
Then, for any $\kappa\in 1\slash 2\,\mathbb{Z}$, cusp forms of weight-$\tilde{\kappa}=(\kappa,\ldots,\kappa)$ 
and nebentypus $\nu^{2\kappa}$ with respect to $\Gamma$ are global sections of the line bundle 
$\omega_{X,\,\nu}^{\otimes 2 \kappa}$. 

\vspace{0.2cm}
Furthermore, for any $f\in\omega_{X,\,\nu}$, i.e., $f$ is a cusp form of weight-$\tilde{1\slash 2}$ and nebentypus $\nu$ with 
respect to $\Gamma$ and $z=(z_1=x_{1}+y_{1},\ldots,z_n=x_{n}+y_{n})$, the Petersson metric on the line 
bundle $\omega_{X,\,\nu}$ is given by 
\begin{align}\label{peterssonip1}
\|f(z)\|_{\omega_{X,\,\nu}}^{2}:=\bigg(\prod_{i=1}^{n}y_{i}^{1\slash 2}\bigg)|f(z)|^{2}. 
\end{align}
Similarly, we assume that $\Omega_{X,\,\nu}$, the line bundle 
of cusp forms of weight-$\tilde{2}$ with nebentypus $\nu$ over $X$ exists. Then, for any $k\in2\mathbb{Z}$, 
cusp forms of weight-$\tilde{k}$ and nebentypus $\nu^{k\slash 2}$ with respect to $\Gamma$ are global sections of the 
line bundle $\Omega_{X,\,\nu}^{\otimes k\slash 2}$. 

\vspace{0.2cm}
Furthermore, for any $f\in\Omega_{X,\,\nu}$, i.e., $f$ is a cusp form of weight-$\tilde{2}$ and nebentypus $\nu$ 
with respect to $\Gamma$ and $z=(z_1=x_{1}+y_{1},\ldots,z_n=x_{n}+y_{n})$, the Petersson metric on the 
line bundle $\Omega_{X,\,\nu}$ is given by 
\begin{align*}
 \|f(z)\|_{\Omega_{X,\,\nu}}^{2}:=\bigg(\prod_{i=1}^{n}y_{i}^{2}\bigg)|f(z)|^{2}. 
\end{align*}
\begin{rem}
For any $z\in X$ and $\kappa\in 1\slash2\,\mathbb{Z}$, from the definition of the Bergman kernel 
$\mathcal{B}_{X,\Omega_{X,\,\nu}}^{2\kappa}(z)$ for the line bundle $\omega_{X, \,\nu}^{\otimes 2\kappa}$ 
from equation \eqref{bkdefn}, we have
\begin{align}\label{remeqn2}
\mathcal{B}_{X,\omega_{X,\,\nu}}^{2\kappa}(z)=\bkappav(z). 
\end{align}
Similarly, for any $z\in X$ and $k\in 2\mathbb{Z}$, from the definition of the Bergman kernel $\mathcal{B}_{X,\Omega_{X,
\,\nu}}^{k\slash 2}(z)$ for the line bundle $\Omega_{X,\,\nu}^{\otimes k\slash 2}$ from 
equation \eqref{bkdefn}, we have
\begin{align*}
\mathcal{B}_{X,\Omega_{X,\,\nu}}^{k\slash 2}(z)=\bkv(z). 
\end{align*}
\end{rem}

\vspace{0.15cm}
\begin{thm}\label{thm1}
Let notation be as above, and we assume that $\Gamma$ is of Type (1), i.e., $X$ is compact. 
Then, for $\kappa\in 1\slash 2\,\mathbb{Z}$, we have
\begin{align*}
\lim_{\kappa\rightarrow \infty}\sup_{z\in X}\frac{1}{\kappa^{n}}\bkappav(z)=O(1);
\end{align*}
and for $k\in 2\mathbb{Z}$, we have
\begin{align*}
 \lim_{k\rightarrow \infty}\sup_{z\in X}\frac{1}{k^{n}}\bkv(z)=O(1),
\end{align*}
where the implied constants in both the above equations are independent of $\Gamma$.   
\begin{proof}
We prove the theorem for $\kappa\in 1\slash 2\,\mathbb{Z}$, and the case for $k\in2\mathbb{Z}$ follows 
automatically with notational changes. For any $z\in X$, from the definition of the curvature form from equation 
\eqref{curvatureform}, and from the definition of the Petersson inner-product for the line 
bundle $\omega_{X,\nu}$ from equation \eqref{peterssonip1}, we have
\begin{align}\label{proofeqn1}
 c_{1}(\omega_{X,\,\nu})(z)=-\frac{i}{2\pi}\partial\overline{\partial}\log\bigg(\prod_{i=1}^{n}y_{i}^{2}
 \bigg)=-\frac{i}{2\pi}\sum_{i=1}^{n}\partial\overline{\partial}\log\big(y_{i}^{2}\big).
\end{align}
For any $1\leq i\leq n$, we compute 
\begin{align}\label{proofeqn2}
\frac{i}{2\pi}\partial\overline{\partial}\log\big(y_{i}^{1\slash 2}\big) =
\frac{i}{4\pi}\partial\overline{\partial}\log\left(\frac{z_{i}-\overline{z}_{i}}{2i}\right)=
-\frac{i}{4\pi}\partial\left(\frac{d\overline{z}_{i}}{z_{i}-\overline{z}_{i}}\right)\notag\\=
\frac{i}{4\pi}\cdot\frac{dz_{i}\wedge d\overline{z}_{i}}{\big(z_{i}-\overline{z}_{i}\big)^{2}}
=-\frac{i}{16\pi}\cdot\frac{dz_{i}\wedge d\overline{z}_{i}}{ y_{i}^{2}}=-\frac{1}{8\pi}\hyp(z_{i}).
\end{align}
Combining equations \eqref{proofeqn1} and \eqref{proofeqn2}, we arrive at
\begin{align*}
 c_{1}(\omega_{X,\,\nu})(z)=\frac{1}{8\pi}\sum_{i=}^{n}\hyp(z_{i})=\frac{1}{8\pi}\hypn(z),
\end{align*}
which shows that the line bundle $\omega_{X,\,\nu}$ is positive, and 
\begin{align}\label{proofeqn3}
\mathrm{det}_{\hypn}\big(c_{1}(
\omega_{X,\,\nu})(z)\big)=(1\slash 8\pi)^{n}. 
\end{align}
So applying estimate \eqref{boucheeqn2} to the complex manifold $X$ with its natural Hermitian metric 
$\hypn$ and the line bundle $\omega_{X,\,\nu}^{\otimes 2\kappa}$, and using equations 
\eqref{remeqn2} and \eqref{proofeqn3}, we find that
\begin{align*}
\lim_{\kappa\rightarrow \infty}\frac{1}{\kappa^{n}}\bkappav(z)=\lim_{\kappa\rightarrow \infty}\frac{1}{
\kappa^{n}}\mathcal{B}_{X,\omega_{X,\,\nu}}^{2\kappa}(z)=O
\bigg(2\mathrm{det}_{\hypn}\big(c_{1}(\omega_{X,\,\nu})(z)\big)\bigg)=O(1).
\end{align*}
As the above limit convergences uniformly in $z\in X$, and as $X$ is compact, we have
\begin{align*}
\sup_{z\in X} \lim_{\kappa\rightarrow \infty}\frac{1}{\kappa^{n}}\bkappav(z)=\lim_{\kappa\rightarrow \infty}
\sup_{z\in X}\frac{1}{\kappa^{n}}\bkappav(z)=O(1),
\end{align*}
which completes the proof of the theorem.
\end{proof}
\end{thm}

\vspace{0.15cm}
\begin{cor}\label{cor1}
Let notation be as above, and we assume that $\Gamma$ is of Type (2), i.e., $X$ is a noncompact complex 
manifold with cusps. Then, for $z\in X$ and $\kappa\in 1\slash 2\,\mathbb{Z}$, 
we have 
\begin{align*}
\limsup_{\kappa\rightarrow\infty}\frac{1}{\kappa^{n}}\bkappav(z)\leq\frac{1}{(8\pi)^{n}};
\end{align*}
and for  $k\in 2\mathbb{Z}$, we have 
\begin{align*}
\limsup_{k\rightarrow \infty}\frac{1}{k^{n}}\bkv(z)\leq\frac{1}{(2\pi)^{n}}.
\end{align*}
\begin{proof}
The proof of the theorem follows from estimate \eqref{bermaneqn}, and from similar 
arguments as in Theorem \ref{thm1}. 
\end{proof}
\end{cor}

\vspace{0.15cm}
\begin{cor}
Let notation be as above, and we assume that $\Gamma$ is of Type (2), i.e., $X$ is a noncompact complex 
manifold with cusps. Let $A$ be a compact subset of $X$. Then, for any $z\in A$ and $\kappa\in 1\slash 2\,
\mathbb{Z}$, we have the following estimate
\begin{align*}
\lim_{\kappa\rightarrow\infty}\frac{1}{\kappa^{n}}\bkv(z)= O_{A}(1);
\end{align*}
and for  $k\in 2\mathbb{Z}$, we have 
\begin{align*}
\lim_{k\rightarrow\infty}\frac{1}{k^{n}}\bkv(z)= O_{A}(1);
\end{align*}
and the implied constant depends on $A$.
\begin{proof}
We prove the theorem for $\kappa\in 1\slash 2\,\mathbb{Z}$, and the case for $k\in2\mathbb{Z}$ follows 
automatically with notational changes. From Corollary \ref{cor1}, for any $z\in A$ and $\kappa\in 
1\slash 2\,\mathbb{Z}$ , we have
\begin{align*}
\limsup_{\kappa\rightarrow\infty}\frac{1}{\kappa^{n}}\bkappav(z)\leq\frac{1}{(8\pi)^{n}}. 
\end{align*}
As $A\subset X$ is compact, we can find a constant $C$ (which depends on $A$) such that 
\begin{align*}
\lim_{\kappa\rightarrow\infty}\frac{1}{\kappa^{n}}\bkappav(z)\leq C, 
\end{align*}
which completes the proof of the corollary.  
\end{proof}
\end{cor}
\paragraph{Acknowledgements}
The author would like to thank J.~Kramer and J.~Jorgenson for introducing him to the area of automorphic forms 
and heat kernels. The author would like to express his gratitude to Archana S. Morye, for many helpful 
discussions and remarks.
{\small{}}

\vspace{0.3cm}
{\small{Anilatmaja Aryasomayajula\\ {\it{anilatmaja@gmail.com}}  \\ 
Department of Mathematics, \hfill\\University of Hyderabad, \\Prof. C.~R.~Rao Road, Gachibowli,\\
Hyderabad, 500046, India}}

\begin{thebibliography}{99}
\addcontentsline{toc}{chapter}{Bibliography}
\bibitem{abbes}  A. Abbes and E. Ullmo, Comparison des metriques d'Arakelov et de Poincare sur $X_0(N)$, Duke Math. J. 80 
(1995), 295--307. 
\bibitem{anil} A. Aryasomayajula, Heat kernel approach for sup-norm bounds for cusp forms of integral and half-integral weight, 
arxiv preprint 
\bibitem{berman} R. J. Berman, Bergman kernels and local holomorphic Morse inequalities, Math. Z. 248 (2004), 
325--344.
\bibitem{bouche} T.~Bouche, Asymptotic results for Hermitian line bundles over complex manifolds: The heat kernel 
approach, Higher-dimensional complex varieties, 67--81, de Gruyter, Berlin, 1996.
\bibitem{codgell} J. W.~Codgell and W.~Luo, The Bergman kernel and M\"ass Equidistribution on the Siegel modular Variety, $\mathrm{Sp}_{2n}
(\mathbb{Z})\backslash\mathfrak{H}_{n}$, Forum Mathematicum 23 (2011), 141--159.
\bibitem{temp} G. Harcos and N. Templier, On the sup-norm of Maass cusp forms of large level. III, 
Math. Ann. 356 (2013), 209--216.
\bibitem{kiral} E. M. Kiral, Bounds on Sup-norms of Half Integral Weight Modular Forms, Acta Arithmetica 165 (2014), 
385-399.
\bibitem{jk2} J.~Friedman, J.~Jorgenson, and J.~Kramer, Uniform sup-norm bounds on average for
cusp forms of higher weights, arXiv preprint arXiv:1305.1348 (2013).
\bibitem{freitag} E. Freitag, Hilbert Modular Forms, Springer-Verlag, Berlin, 1990. 
\bibitem{jk1} J.~Jorgenson and J.~Kramer, Bounding the sup-norm of automorphic
forms, GAFA 14 (2004), 1267--1277.
\bibitem{liu} S-C. Liu, Equidistribution of Hecke eigenforms on the Hilbert modular varieties, Journal of 
Number Theory 127 (2007),  1--9.
\bibitem{raphael} R. S. Steiner, Uniform bounds on sup-norms of holomorphic forms of real weight, arxiv preprint 
arXiv:1406.2918.
\end{thebibliography}
\end{document}